\documentclass{article}\begin{document}
\centerline{\bf The complete evaluation of Rogers Ramanujan and other}
\centerline{\bf continued fractions with elliptic functions}
\vskip .4in

\centerline{Nikos Bagis}
\centerline{Department of Informatics, Aristotle University}
\centerline{Thessaloniki, Greece}
\centerline{bagkis@hotmail.com}
\vskip .2in

\[
\]
\textbf{Keywords}: Ramanujan; Continued Fractions; Elliptic Functions; Modular Forms
\[
\]
\centerline{\bf Abstract}

\begin{quote}
In this article we present evaluations of continued fractions studied by Ramanujan. More precisely we give the complete polynomial equations of Rogers-Ramanujan and other continued fractions, using tools from the elementary theory of the Elliptic functions. We see that all these fractions are roots of polynomials with coeficients depending only on the inverse elliptic nome-q and in some cases the Elliptic Integral-K.
In most of simplifications of formulas we use Mathematica.   

\end{quote}

\newpage

\section{Introductory Definitions, and some Properties of the Rogers Ramanujan Continued Fraction}

For $\left|q\right|<1$, the Rogers Ramanujan continued fraction (RRCF) is defined as
\begin{equation}
R(q):=\frac{q^{1/5}}{1+}\frac{q^1}{1+}\frac{q^2}{1+}\frac{q^3}{1+}\cdots  
\end{equation}
We also define
\begin{equation}
(a;q)_n:=\prod^{n-1}_{k=0}{(1-aq^k)}
\end{equation}
\begin{equation}
f(-q):=\prod^{\infty}_{n=1}(1-q^n)=(q;q)_{\infty}
\end{equation}
\begin{equation}
\Phi(-q):=\prod^{\infty}_{n=1}(1+q^n)=(-q;q)_{\infty}
\end{equation}
and also hold the following relations of Ramanujan
\begin{equation}
\frac{1}{R(q)}-1-R(q)=\frac{f(-q^{1/5})}{q^{1/5}f(-q^5)} 
\end{equation} 
\begin{equation}
\frac{1}{R^5(q)}-11-R^5(q)=\frac{f^6(-q)}{q f^6(-q^5)} 
\end{equation} 
\[
\]
From the Theory of Elliptic Functions we have:\\
Let
\begin{equation}
K(x)=\int^{\pi/2}_{0}\frac{1}{\sqrt{1-x^2\sin(t)^2}}dt
\end{equation}
It is known that the inverse elliptic nome $k=k_r$, $k'^2_r=1-k^2_r$ is the solution of
\begin{equation}
\frac{K\left(k'_r\right)}{K(k)}=\sqrt{r}
\end{equation}
We also denote $K[r]=K(k_r)$. In what it follows we assume that $r \in \bf R^{*}_+ \rm$. When $r$ is rational then $k_r$ is algebraic.   
\begin{equation}
k_r=\frac{8q^{1/2}\Phi(-q)^{12}}{1+\sqrt{1+64q\Phi(-q)^{24}}}
\end{equation}
We can write the functions $f$ and $\Phi$ using elliptic functions. It holds
\begin{equation}
\Phi(-q)=\frac{2^{5/6}q^{-1/24}(k_r)^{1/12}}{(k'_r)^{1/6}}
\end{equation}
\begin{equation}
f(-q)^8=\frac{2^{8/3}}{\pi^4}q^{-1/3}(k_r)^{2/3}(k'_r)^{8/3}K(k_r)^4
\end{equation}
also holds 
\begin{equation}
f(-q^2)^6=\frac{2k_r k'_r K(k_r)^3}{\pi^3 q^{1/2}}
\end{equation}
From [19] it is known that
\begin{equation}
R'(q)=1/5q^{-5/6}f(-q)^4R(q)\sqrt[6]{R(q)^{-5}-11-R(q)^5}
\end{equation}

\section{Evaluations of Rogers Ramanujan Continued Fraction}

\textbf{Theorem 2.1}\\
If $q=e^{-\pi\sqrt{r}}$ and 
\begin{equation}
a=a_r=\left(\frac{k'_r}{k'_{25r}}\right)^2\sqrt{\frac{k_r}{k_{25r}}}M_5(r)^{-3}
\end{equation}
Then
\begin{equation}
R(q)=\left(-\frac{11}{2}-\frac{a_r}{2}+\frac{1}{2}\sqrt{125+22a_r+a^2_r}\right)^{1/5} 
\end{equation}
Where $M_5(r)$ is root of: $(5x-1)^5(1-x)=256 (k_r)^2 (k'_r)^2 x$.\\
\textbf{Proof.}\\
Suppose that $N=n^2\mu$, where $n$ is positive integer and $\mu$ is positive real then it holds that
\begin{equation}
K[n^2\mu]=M_n(\mu)K[\mu]
\end{equation}
Where $K[\mu]=K(k_{\mu})$\\ 
The following formula for $M_5(r)$ is known
\begin{equation}
(5M_5(r)-1)^5(1-M_5(r))=256(k_r)^2 (k'_r)^2M_5(r)
\end{equation}
Thus if we use (6) and (11) and the above consequence of the Theory of Elliptic Functions, we get:
\[
R^{-5}(q)-11-R^{5}(q)=\frac{f^6(-q)}{q f^6(-q^5)}=a=a_r
\]
Solving with respect to $R(q)$ we get the result.
\[
\]
The relation between $k_{25r}$ and $k_r$ is
\begin{equation}
k_rk_{25r}+k'_rk'_{25r}+2\cdot4^{1/3} (k_rk_{25r}k'_rk'_{25r})^{1/3}=1
\end{equation}
\[
\]
We will try to evaluate $k_{25r}$. For this we set 
\begin{equation}
k_{25r}k_r=w^2
\end{equation} 
Then setting directly to (18) the following parametrization of $w$ (see also [8] pg.280):
\begin{equation} 
w=\sqrt{\frac{L(18+L)}{6(64+3L)}}
\end{equation} 
we get
\begin{equation} 
\frac{(k_{25r})^{1/2}}{w^{1/2}}=\frac{w^{1/2}}{(k_r)^{1/2}}=\frac{1}{2}\sqrt{4+\frac{2}{3}\left(\frac{L^{1/6}}{M^{1/6}}-4\frac{M^{1/6}}{L^{1/6}}\right)^2}+\frac{1}{2}\sqrt{\frac{2}{3}}\left(\frac{L^{1/6}}{M^{1/6}}-4\frac{M^{1/6}}{L^{1/6}}\right)
\end{equation} 
where $$M=\frac{18+L}{64+3L}$$
From the above relations we get also
\begin{equation}
-\frac{k_r-w}{\sqrt{k_rw}}=\frac{k_{25r}-w}{\sqrt{k_{25r}w}}=\sqrt{\frac{2}{3}}\left(\frac{L^{1/6}}{M^{1/6}}-4\frac{M^{1/6}}{L^{1/6}}\right)
\end{equation}
We can consider the above equations as follows: Taking an arbitrary number $L$ we construct an $w$. Now for this $w$ we evaluate the two numbers $k_{25r}$ and $k_r$. Thus when we know the $w$, the $k_{r}$ and $k_{25r}$ are given by (20) and (21). The result is: We can set a number $L$ and from this calculate the two inverse elliptic nome`s, or equivalently, find easy solutions of (18). But we don't know the $r$. One can see (from the definition of $k_r$) that the $r$ can be evaluated from equation 
\begin{equation}
r=r_{k_1}=r[L]=\frac{K^2(k'_r)}{K^2(k_r)}
\end{equation}
or
\begin{equation}
r=r_{k_{25}}=r[25L]=\frac{1}{25}\frac{K^2(k'_{25r})}{K^2(k_{25r})}
\end{equation}
However there is no way to evaluate the $r$ in a closed form, such as roots of polynomials, or else. Some numerical evaluations as we will see, indicate as that even $k_r$ are algebraic numbers, the $r$ are not rational.\\  
\textbf{Theorem 2.2}\\
Set
\begin{equation}
A_L=a\left(\frac{K^2(k'_L)}{K^2(k_L)}\right)=\frac{(k_L)^3 (1-(k_{L})^2) M_5(L)^{-3}}{(k_L)^2w-w^5}
\end{equation}
then
\begin{equation}
R\left(e^{-\pi \sqrt{r[L]}}\right)=\left(-\frac{11}{2}-\frac{A_L}{2}+\frac{1}{2}\sqrt{125+22A_L+A^2_L}\right)^{1/5} 
\end{equation}
where the $k_L$ and $w$ are given by (19) and (20).\\
\textbf{Example.}\\
Set $L=1/3$ then $$w=\frac{1}{3}\sqrt{\frac{11}{78}}$$ and $$k_{25r}=\frac{1}{3}\sqrt{\frac{11}{78}}\left(\frac{-4(\frac{11}{13})^{1/6}+(\frac{13}{11})^{1/6}}{\sqrt{6}}+\frac{1}{2}\sqrt{4+\frac{2}{3}\left(-4\left(\frac{11}{13}\right)^{1/6}+\left(\frac{13}{11}\right)^{1/6}\right)^2}\right)^2$$
and
$$k_r=\frac{\frac{1}{3}\sqrt{\frac{11}{78}}}{\left(\frac{-4(\frac{11}{13})^{1/6}+(\frac{13}{11})^{1/6}}{\sqrt{6}}+\frac{1}{2}\sqrt{4+\frac{2}{3}\left(-4\left(\frac{11}{13}\right)^{1/6}+\left(\frac{13}{11}\right)^{1/6}\right)^2}\right)^2}$$
where the $r$ is given by
$$r=\frac{K^2\left(\sqrt{1-k_r^2}\right)}{K^2\left(k_r\right)}$$
\[
\]
Now we can see (The results are known in the Theory of Elliptic Functions) how we can found evaluations of $R(q)$ when $r$ is given and $k_r$ is known:\\
From (20) it is
\begin{equation} 
L=-9+9w^2+\sqrt{3}\sqrt{27+74w^2+27w^4} 
\end{equation} 
from the relation between $M$ and $L$ we get
\begin{equation}  
M=\frac{1}{64}\left(9-9w^2+\sqrt{81+222w^2+81w^4}\right)
\end{equation}
Hence from (22)
\begin{equation}
t=\frac{w-k_r}{\sqrt{k_rw}}
\end{equation} 
also
\begin{equation}
t=\sqrt{\frac{2}{3}}\left(\frac{1}{y^{1/6}}-4y^{1/6}\right)
\end{equation}
where $y=M/L$. Hence ($k=k_r$):
\begin{equation} 
\frac{M}{L}=\left(\frac{\sqrt{3}(k-w)+\sqrt{3k^2+26kw+3w^2}}{8\sqrt{2kw}}\right)^6
\end{equation}
or

$$\frac{1}{64}\left(\frac{9-9w^2+\sqrt{81+222w^2+81w^4}}{-9+9w^2+\sqrt{81+222w^2+81w^4}}\right)
=\left(\frac{\sqrt{3}(k-w)+\sqrt{3k^2+26kw+3w^2}}{8\sqrt{2kw}}\right)^6$$
or
\begin{equation}
\frac{\sqrt{6}w}{-9+9w^2+\sqrt{81+222w^2+81w^4}}=\left(\frac{\sqrt{3}(k-w)+\sqrt{3k^2+26kw+3w^2}}{8\sqrt{2kw}}\right)^3
\end{equation}
setting now 
\begin{equation}
k^{*}_r=\left(\frac{-1+4p^2+\sqrt{1-2p^2+16p^4}}{\sqrt{6}p}\right)^2
\end{equation}
and
\begin{equation}
w=\left(\frac{6^{3/4}p^{3/2}}{\sqrt{-1+64p^6+\sqrt{1+88p^6+4096p^{12}}}}\right)^2
\end{equation}
$$W=-1+4p^2+\sqrt{1-2p^2+16p^4}$$ 
and
$$T=-1+64p^6+\sqrt{1+88p^6+4096p^{12}}$$
we have 
$$k^{*}_rw=k_r$$ 
Also
\begin{equation}
p=\left(\frac{T(2+T)}{216+128T}\right)^{1/6}=\left(\frac{W(2+W)}{6+8W}\right)^{1/2}
\end{equation}
But 
$$w=6k_r\left(\frac{W+2}{(6+8W)W}\right)\eqno{(35a)}$$
$$T=\frac{\sqrt{6}W^2}{k_r}\sqrt{\frac{W(2+W)}{6+8W}}$$
where the equation for finding $W$ from $k_r$ is 
$$
-108k^2_r\left(\frac{W(W+2)}{8W+6}\right)^{5/2}+\sqrt{6}k_rW^2\left(1-64\left(\frac{W(W+2)}{8W+6}\right)^3\right)+$$
\begin{equation}
+3W^4\left(\frac{W(W+2)}{8W+6}\right)^{1/2}=0
\end{equation} 
We give the complete polynomial equation of $p$ arising from (35):
$$
k^2_r+2\sqrt{6}k_rk'^2_rp-24k^2_rp^2-10\sqrt{6}k_rk'^2_rp^3+
240k^2_rp^4+32\sqrt{6}k_rk'^2_rp^5+$$
$$
+(54-1388k^2_r+54k^4_r)p^6-128\sqrt{6}k_rk'^2_rp^7+3840k^2_rp^8+640\sqrt{6}k_rk'^2_rp^9-$$
\begin{equation}
-6144k^2_rp^{10}-2048\sqrt{6}k_rk'^2_rp^{11}+4096k^2_rp^{12}=0
\end{equation}
It is evident that the Rogers Ramanujan Continued Fraction is  root of polynomial equation with coeficients depending by $k_r$.

From (33) we have $$\sqrt{6k^{*}_r}=\frac{-1+4p^2+\sqrt{1-2p^2+16p^4}}{p}
$$
where $p$ is root of (37).\\ Using Mathematica we get the following simplification formula for
$$x=\sqrt{k^{*}_r}=1/\sqrt{k_{25r}}$$
$$k^2_r+4k'^2_rk_rx-6k^2_rx^2+20k'^2_rk_rx^3+15x^4
-16k'^2_rk_rx^5+(16-52k^2_r+16k^4_r)x^6+$$
\begin{equation}
+16k'^2_rk_rx^7+15k^2_rx^8-20k'^2_rk_rx^9-6k^2_rx^{10}-4k'^2_rk_rx^{11}+k^2_rx^{12}=0
\end{equation}
set now 
\begin{equation}
c_r=\frac{k'^2_r (k^{*}_r)^5}{(k^{*}_r)^4-k^2_r}
\end{equation}
and $$G(q)=(R^{-5}(q)-11-R^5(q))^{1/3}$$
then\\
\textbf{Theorem 2.3}\\
\textbf{i)} 
$$3125c^2_r-6250c^{5/3}_rG(q)+4375c^{4/3}_rG^2(q)-1500c_rG^3(q)+275c^{2/3}_rG^4(q)+$$
$$
+2c^{1/3}_r(-13+128k'^2_rk^2_r)G^5(q)+G^6(q)=0\eqno{(39a)}
$$
Also\\
\textbf{ii)}
$$k^6_r+k^3_r(-16+10k^2_r)w+15k^4_rw^2-20k^3_rw^3+15k^2_rw^4+k_r(10-16k^2_r)w^5+w^6=0\eqno{(39b)}$$
Once we know $k_r$ we can calculate $w$ from the above equation and hence the $k_{25r}$.\\ Hence the problem reduces to solve 6-th degree equations. The first is (17) the second is (39b) and for (RRCF) (39a).\\    
\textbf{Proof.}\\
\textbf{i)} We have
$$R(q)^{-5}-11-R^5(q)=a_r=\frac{k_r^3(1-k_r^2)}{w(k_r^2-w^4)}M_5(r)^{-3}$$
$$=\frac{k'^2_r (k^{*}_r)^5}{(k^{*}_r)^4-k^2_r}M_5(r)^{-3}$$
and $M_5(r)$ satisfies $(5x-1)^5(1-x)=256(k_r)^2(k'_r)^2x$.\\
After elementary algebraic calculations we get the result.\\ 
\textbf{ii)}\\ 
From (35) we get:
$$-\sqrt{6}t-3U+108t^2U^5+64\sqrt{6}tU^6=0$$
and
$$t=\frac{k_r}{W^2}=\frac{w}{6U^2}$$
and
$$-32\sqrt{6}wU^6+(9-54w^2)U^3+3\sqrt{6}w=0\eqno{(a)}$$ 
$$U=\left(\frac{W(W+2)}{8W+6}\right)^{1/2}\eqno{(b)}$$
and also
$$\left(\frac{W(W+2)}{8W+6}\right)^{1/2}=\sqrt{\frac{w}{6k_r}}W\eqno{(c)}$$
Hence solving the system we obtain the 6-th degree equation (39b).\\
\textbf{Corollary 2.1}\\
The solution of (39b) with respect to $k_r$ when we know $w$ is
\begin{equation} 
\frac{w^{1/2}}{(k_r)^{1/2}}=\frac{1}{2}\sqrt{4+\frac{2}{3}\left(\frac{L^{1/6}}{M^{1/6}}-4\frac{M^{1/6}}{L^{1/6}}\right)^2}+\frac{1}{2}\sqrt{\frac{2}{3}}\left(\frac{L^{1/6}}{M^{1/6}}-4\frac{M^{1/6}}{L^{1/6}}\right)
\end{equation}
Where
$$
w=\sqrt{\frac{L(18+L)}{6(64+3L)}}
$$ 
$$M=\frac{18+L}{64+3L}$$
\textbf{Theorem 2.4}
$$R'(q)=\frac{2^{4/3}(k_r)^{5/12}(k'_r)^{5/3}}{5(k_{25r})^{1/12}(k'_{25r})^{1/3}\sqrt{M_5(r)}}\times$$ 
\begin{equation}
\times\left(-\frac{11}{2}-\frac{a_r}{2}+\frac{1}{2}\sqrt{125+22a_r+a^2_r}\right)^{1/5}\frac{K^2(k_r)}{\pi^2 q}
\end{equation}
\textbf{Proof.}\\
Combining (11) and (13) and Theorem 2.1 we get the proof.\\
\textbf{Evaluations.}
\[
R(e^{-2\pi})=\frac{-1}{2}-\frac{\sqrt{5}}{2}+\sqrt{\frac{5+\sqrt{5}}{2}}
\]
\[
R'(e^{-2\pi})=8\sqrt{\frac{2}{5}\left(9+5\sqrt{5}-2\sqrt{50+22\sqrt{5}}\right)} \frac{e^{2\pi}}{\pi^3}\Gamma\left(\frac{5}{4}\right)^4
\]
\[
\]
Sumarizing our results we can say that:\\
1) Theorem 2.1 is quite usefull for evaluating $R(q)$ when we know $k_r$ and $k_{25r}$. But this it whas known allrady to Ramanujan by using the function $X(-q)=(-q;q^2)_{\infty}$, (see [8]).\\ 
2) Theorem 2.2 is more kind of a Lemma rather a Theorem and it might help for further research.\\
3) Theorem 2.3 is a proof that the Rogers Ramanujan continued fraction is a root of a polynomial equation with coeficients the $k_r$ where $r$ positive real. The polynomial equation (39a), is a version of icosahedral equation for the (RRCF) (see[20]).
\\4) Theorem 2.4 is a consequence of a Ramanujan integral first proved by Andrews (see [5]) and it is usefull for evaluations of $R'(e^{-\pi \sqrt{r}})$, $r\in \bf Q\rm$.\\
The above theorems can used to derive also modular equations of $R(q)$, from the modular equations of $k_r$. More precisely we can guess an equality with the help of a methematical pacage (for example in Mathematica there exist the command 'recognize') and then proceed to proof using Theorems which we presented in this article. We follow this prosedure with other continued fractions (the Rogers Ramanujan is litle dificult), such as the cubic or Ramanujan-Gollnitz-Gordon.\\These two last continued fractions are more easy to handle. The elliptic function theory and the sigular moduli $k_r$ will exctrac and give us several proofs of modular indenties.      

\section{The H-Continued Fraction}

Heng Huat Chan and Sen-Shan Huang [11] studied the Ramanujan Gollnitz-Gordon continued fraction
\begin{equation}
H(q):=\frac{q^{1/2}}{1+q+}\frac{q^2}{1+q^3+}\frac{q^4}{1+q^5+}\frac{q^6}{1+q^7+}\ldots
\end{equation} 
where $\left|q\right|<1$.\\
In the paper of C. Adiga and T. Kim [2] one can find the next identity for this fraction   
\begin{equation}
H(q)^{-1}-H(q)=\frac{M^2(q^2)}{M^2(q^4)}
\end{equation}
where
\begin{equation}
M(q)=\frac{q^{1/8}}{1+}\frac{-q^1}{1+q^1}\frac{-q^2}{1+q^2+}\frac{-q^3}{1+q^3+}\ldots=q^{1/8}\frac{(q^2;q^2)_{\infty}}{(q;q^2)_{\infty}}
\end{equation}

Next we will use some properties of the inverse Elliptic Nome and show how this can help us to evaluate the H-fraction. For to complete our purpose we need the relation between $k_r$ and $k_{4r}$. There holds the following\\
\textbf{Lemma 3.1}\\
If $r>0$, then
\begin{equation}
k_{4r}=\frac{1-k'_{r}}{1+k'_{r}}      
\end{equation} 
and 
\begin{equation}
K[4r]=\frac{1+k'_r}{2}K[r]
\end{equation}
\textbf{Proof.}\\
For (45) see ([8], pg. 102, 215). The identity for $K[4r]$ is known from the theory of elliptic functions\\ 
\textbf{Theorem 3.1}
\begin{equation}
H(q)=-P+\sqrt{P^2+1} 
\end{equation} 
where
$$P=\frac{k_r}{(1-k'_r)}$$   
or equivalently 
\begin{equation}
k_r=\frac{4(H-H^3)}{(1+H^2)^2}
\end{equation}
\textbf{Proof.}\\
It is known that, under some conditions in the sequence $b_n$ (see [16]) it holds
\begin{equation}
\frac{1}{1+}\frac{-b_1}{1+b_1}\frac{-b_2}{1+b_2+}\ldots=1+\sum^{\infty}_{n=1}\prod^{n}_{k=1}b_k
\end{equation}
Hence if we set $b_n=q^n$, $\left|q\right|<1$, then 
\begin{equation}
M(q)=\theta_2(q^{1/2})=q^{-1/8}\sqrt{\frac{k_{r/4}K(k_{r/4})}{2\pi}}
\end{equation}
where $$\theta_2(q)=\sum^{\infty}_{n=-\infty}q^{(n+1/2)^2}$$ 
$$\sum^{\infty}_{n=0}q^{n(n+1)/2}=1/2q^{-1/8}\theta_2(q^{1/2})$$
Using Lemma 3.1 and identity (43) we get the proof.\\
\textbf{Theorem 3.2}\\
If $ab=\pi^2$, then
\begin{equation}\
\left(H(e^{-a})+2-\frac{1}{H(e^{-a})}\right)\left(H(e^{-b})+2-\frac{1}{H(e^{-b})}\right)=8
\end{equation}
\textbf{Proof.}\\
Set 
\begin{equation}
\psi(q)=\sum^{\infty}_{n=0}q^{(n+1)n/2}
\end{equation}
and
\begin{equation}
\phi(q)=\sum^{\infty}_{n=-\infty}q^{n^2}
\end{equation}
Identity (51) becomes.\\
If $ab=4\pi^2$
\begin{equation}
\left(2-\frac{\psi(e^{-a})^2}{e^{-a/4}\psi(e^{-2a})^2}\right)\left(2-\frac{\psi(e^{-b})^2}{e^{-b/4}\psi(e^{-2b})^2}\right)=8
\end{equation}
From [8] pg.43 we have if $ab=2\pi$, then
\begin{equation}
\psi(e^{-a^2})=\frac{\sqrt{b}}{2\sqrt{a}}e^{a^2/8}\phi(-e^{-b^2/2})
\end{equation} 
\begin{equation}
\psi(e^{-2a^2})=\frac{\sqrt{b/2}}{2\sqrt{a}}e^{a^2/4}\phi(-e^{-b^2/4})
\end{equation}
Hence if $ab=\pi^2/4$, ([8] pg.98)
\begin{equation}
\left(1-\frac{\phi(e^{-a})}{\phi(-e^{-a})}\right)\left(1-\frac{\phi(e^{-b})}{\phi(-e^{-b})}\right)=2
\end{equation}
But this is equivalent to 
\begin{equation}
k'_{1/(4r)}=\frac{1-k'_r}{1+k'_r}
\end{equation}
(For details [8] pg.98, 102 and 215),
which is equivalent to
\begin{equation}
k_r=k'_{1/r}.
\end{equation} 
But this is true from the definition of the modulus-$k$ (see relation (8)). 
This completes the proof.\\
\textbf{Corollary 3.1}\\
If $ab=\pi^2$
\begin{equation} 
(1+\sqrt{2}+H(e^{-a}))(1+\sqrt{2}+H(e^{-b}))=2(2+\sqrt{2})
\end{equation}
\textbf{Proof.}\\
This follows from Theorem 3.2 and as in [8] pg. 84\\ 
\textbf{Evaluations.}
\[
H\left(e^{-\pi/2}\right)=\sqrt{1+2\sqrt{2}-2\sqrt{2+\sqrt{2}}}
\]
\[
H\left(e^{-\pi\frac{\sqrt{2}}{2}}\right)=\sqrt{3+2\sqrt{2}-2\sqrt{4+3\sqrt{2}}}
\]

Now it is easy to see how we can construct modular equations of these continued fractions from the modular equations of the inverse elliptic nome. For example for the $H$ continued fraction we give the second degree modular equation:\\
\textbf{Theorem 3.3}  
\begin{equation}
H^2(q)=\frac{H(q^2)-H^2(q^2)}{1+H(q^2)}
\end{equation}
\textbf{Proof.}\\
If $ab=4$ and $k_r=k(e^{-\pi\sqrt{r}})$, $q_r=e^{-\pi\sqrt{r}}$ then
\[
(1+k_a)(1+k_b)=2
\] 
which can be written as
\[
(1+k_{4a})(1+k'_a)=2
\] 
But from Theorem 3.1
\[
k_a=\frac{4H^2(q^{1/2})}{(H^2(q^{1/2})-1)^2}
\]
and the result follows after elementary algebraic computations.
\[
\]
Also from (48) and (61) one can get
\begin{equation} 
\sqrt{k'_r}=\frac{H(q^2)+2H(q^2)-1}{H(q^2)-2H(q^2)-1}
\end{equation}

For to proceed we must mention that the relation between $k_{9r}$ and $k_r$ is given by (see [8]):
\begin{equation}
\sqrt{k_rk_{9r}}+\sqrt{k'_rk'_{9r}}=1
\end{equation}

\section{The Ramanujan's Cubic Continued Fraction}

Let 
\begin{equation} V(q):=\frac{q^{1/3}}{1+}\frac{q+q^2}{1+}\frac{q^2+q^4}{1+}\frac{q^3+q^6}{1+}\ldots
\end{equation}
is the Ramanujan's cubic continued fraction, then\\ 
\textbf{Lemma 4.1}
\begin{equation}
V(q)=\frac{2^{-1/3}(k_{9r})^{1/4}(k'_{r})^{1/6}}{(k_r)^{1/12}(k'_{9r})^{1/2}}
\end{equation}
where the $k_{9r}$ are given by (63).\\
\textbf{Proof.}\\
It is known (see and [14] pg. 596) that
\[
V(q)=q^{1/3}\frac{(q;q^2)_{\infty}}{(q^3;q^6)^3_{\infty}}
\] 
But 
$$\Phi(-q)=(-q,q)_{\infty}=\frac{1}{(q,q^2)_{\infty}}$$
thus
$$V(q)=q^{1/3}\frac{\Phi(-q^3)^3}{\Phi(-q)}$$
and equation (65) follows from (10).\\
\textbf{Lemma 4.2}\\
If $$G(x)=\frac{x}{\sqrt{2\sqrt{x}-3x+2x^{3/2}-2\sqrt{x}\sqrt{1-3\sqrt{x}+4x-3x^{3/2}+x^2}}}$$ 
and $w$ is defined by 
$$k_{9r}=\frac{w}{k_r}$$
then
$$k'_{9r}=\frac{(1-\sqrt{w})^2}{k'_r}$$
and
\begin{equation}
k_r=G(w)
\end{equation}
\textbf{Proof.}\\
Set the values of $k_r$ and $k_{9r}$ in (63).\\
\\

Also holds
$$\frac{1}{(V(q)V(q^3))^{12}}=256w\frac{\left(1-\frac{w^2}{G(w)^2}\right)^2}{1-G(w)^2}\frac{\left(1-\frac{G(w)^2G^{-1}(w/G(w))^2}{w^2}\right)^3}{G^{-1}(w/G(w))^3}$$
If we set
\begin{equation} 
W=2-3\sqrt{w}+2w-2(1-\sqrt{w})\sqrt{1-\sqrt{w}+w}
\end{equation}
then
\begin{equation}
V(q)=\frac{(k'_r)^{2/3}w^{1/4}}{2^{1/3}(k_r)^{1/3}(1-\sqrt{w})}=\frac{(W-w^{3/2})^{1/3}W^{-1/6}}{2^{1/3}(1-\sqrt{w})} 
\end{equation}
after solving (67) with respect to $w$ and making the simplifications we arrive at
\begin{equation}
2V^3(q)=\frac{\sqrt{W}}{(1+\sqrt{W})^2} 
\end{equation}
and
\begin{equation}
(k_r)^2=\sqrt{W}\left(\frac{2+\sqrt{W}}{1+2\sqrt{W}}\right)^3
\end{equation}
Hence we get the following equation 
\begin{equation}
(k_r)^{2/3}=Z^{2}\frac{\sqrt{2}V(q)^{3/2}+Z^3}{-\sqrt{2}V(q)^{3/2}+2Z^3}
\end{equation}
Where $Z=\sqrt[12]{W(q)}$. Reducing the above equation in polynomial form we have
\begin{equation}
sk_r^{2/3}+sZ^2-2k_r^{2/3}Z^3+Z^5=0
\end{equation} 
and
\begin{equation}
s^2=2V^3(q)=\frac{Z^6}{(1+Z^6)^2}
\end{equation}
From these two last equations we arrive to\\
\textbf{Theorem 4.1}\\
Set $T=\sqrt{1-8V(q)^3}$ then holds the next equation
\begin{equation}
(k_r)^2=\frac{(1-T)(3+T)^3}{(1+T)(3-T)^3}
\end{equation}    
\textbf{Corollary 4.1}\\
If $X=\sqrt{W(q)}=\frac{1-T}{1+T}$ and $Y=\sqrt{W(q^2)}$, then 
\begin{equation}
X^{1/2}\left(\frac{2+X}{1+2X}\right)^{3/2}=2\frac{Y^{1/4}}{\left(\frac{1+2Y}{2+Y}\right)^{3/4}+Y^{1/2}\left(\frac{2+Y}{1+2Y}\right)^{3/4}}
\end{equation}

The dublication formula is\\
\textbf{Proposition 4.1}\\
Set $u=T(q^2)$, $v=T(q)$, then
\begin{equation}
\frac{\sqrt{(1-u)}(3+u)^{3/2}}{\sqrt{(1+u)}(3-u)^{3/2}}=\frac{(3-v)^{3/2}\sqrt{1+v}-4v^{3/2}}{(3-v)^{3/2}\sqrt{1+v}+4v^{3/2}} 
\end{equation}
We can simplify the problem of finding modular equations of degree 3 using the Cubic continued fraction. As someone can see with direct algebraic evaluations and with definitions of $W$, $V(q)$ and Lemma 4.2 there holds:\\  
\textbf{Proposition 4.2}\\
If
$$k_{9r}=\frac{w}{k_r}$$
then
\begin{equation}
w=\left(\frac{1-4V(q)^3-8V(q)^6-\sqrt{1-8V(q)^3}}{4V(q)^3\left(1-2V(q)^3-\sqrt{1-8V(q)^3}\right)}\right)^2
\end{equation}
The Ramanujan's modular equation which relates $V(q)$ and $V(q^3)$ is
\begin{equation}
V(q)^3=V(q^3)\frac{1-V(q^3)+V(q^3)^2}{1+2V(q^3)+4V(q^3)^2}
\end{equation}
(see [10]) one can get from the above formula and Proposition 4.2 the folowing:\\
\textbf{Proposition 4.3}
\begin{equation}
k_{81r}=\left(\frac{1+2V(q^3)^2-\sqrt{1-8V(q^3)^3}}{1+2V(q^3)^2+\sqrt{1-8V(q^3)^3}}\right)^2k_r
\end{equation}    
\textbf{Corollary 4.2}\\
Set $u=H(q)$, $v=H(q^6)$ and $$t=\frac{4T(q)}{(1+T(q))(3-T(q))}$$, then 
\begin{equation}
\frac{\sqrt{u^4-6u^2+1}(v^2+2v-1)}{(u^2+1)(v^2-2v-1)}=t
\end{equation}
\textbf{Proof.}\\
Set $q\rightarrow q^3$ in (61) and then use Proposition 4.2 .\\  
\textbf{Evaluations}\\
a)
$$V(e^{-\pi})=\frac{1}{2^{2/3}}\left(-67-39\sqrt{3}+(9+6\sqrt{3})\sqrt{2(12+7\sqrt{3})}\right)$$
b)
$$
\frac{(3-T(e^{-\pi\sqrt{2}}))^3(3+T(e^{-\pi\sqrt{2}}))^3}{(1-T(e^{-\pi\sqrt{2}}))(1+T(e^{-\pi\sqrt{2}}))}=5832
$$
Using the tables of $k_r$ we can find a wide number of evaluations for the cubic continued fraction.\\
c) From [13] we have 
$$b_{M,N}=\frac{Ne^{-\frac{(N-1)\pi}{4}\sqrt{\frac{M}{N}}}\psi^2(-e^{-\pi\sqrt{MN}})\phi^2(-e^{-2\pi\sqrt{MN}})}{\psi^2(-e^{-\pi\sqrt{\frac{M}{N}}})\phi^2(-e^{-2\pi\sqrt{\frac{M}{N}}})}$$
Then from Theorem 4.1 
$$b^2_{M,3}=\frac{9(1-T^2)}{T^2(9-T^2)}$$
where
$$\phi(q):=\sum^{\infty}_{n=-\infty}q^{n^2}$$
$$\psi(q):=\sum^{\infty}_{n=0}q^{n(n+1)/2}$$
$\left|q\right|<1$.

\section{Other Continued Fractions}

\textbf{Section 1.}

Another continued fraction is
\begin{equation}
S(q)=q^{1/8}\frac{1}{1+}\frac{q}{1+}\frac{q^2+q}{1+}\frac{q^3}{1+}\frac{q^4+q^2}{1+}\frac{q^5}{1+}\ldots 
\end{equation}
for which it is known that
\begin{equation}
S(q)=q^{1/8}\frac{(-q^2;q^2)_{\infty}}{(-q;q^2)_{\infty}}
\end{equation}
after using Euler's Theorem: $(-q,q)_{\infty}=1/(q;q^2)_{\infty}$ and making some simplifications and rearrangements in the products we find 
$$S(q)=q^{1/8}\frac{\Phi(-q^2)^2}{\Phi(-q)}$$
Now making use of (10) we get
\begin{equation}
S(q)=\frac{2^{-1/6}(k_{4r})^{1/6}(k'_r)^{1/6}}{(k_r)^{1/12}(k'_{4r})^{1/3}}
\end{equation}
Using the relation between $k_{4r}$ and $k_r$ form Lemma 3.1, 
we get\\
\textbf{Theorem 5.1}
\begin{equation}
S(q)=\frac{(k_r)^{1/4}}{\sqrt{2}}
\end{equation}
Hence the fraction $S$ is the inverse elliptic nome and as someone can see there holds a very large number of modular equations, but since it is trivial we not mention here.\\
\textbf{Section 2.}\\ The continued fraction
\begin{equation}
Q(q)=\frac{q^{1/2}}{1-q+}\frac{q(1-q)^2}{(1-q)(q^2+1)+}\frac{q(1-q^3)^2}{(1-q)(q^4+1)+}\frac{q(1-q^5)^2}{(1-q)(q^6+1)+}\ldots
\end{equation} 
which is known that  
\begin{equation}
Q(q)=q^{1/2}\frac{(q^4;q^4)^2_{\infty}}{(q^2;q^4)^2_{\infty}}=M(q^2)^2
\end{equation}
it becomes
\begin{equation}
Q(q)=...=q^{1/2}f(-q^4)^2\Phi(-q^2)^2=M(q^2)^2
\end{equation}
or\\
\textbf{Theorem 5.2}
\begin{equation}
Q(q)=\frac{1}{\pi}K(k_{4r})\sqrt{k_{4r}}=\frac{1}{2\pi}K(k_r)k_r 
\end{equation}
\textbf{Proof.}\\
It follows from the relation between $Q$ and $M$.\\
\textbf{Evaluation.}
$$Q(e^{-\pi \sqrt{2}})=\frac{\sqrt{2}-1}{\sqrt{2\pi}}\frac{\Gamma(9/8)}{\Gamma(5/8)}$$
\textbf{Theorem 5.3}\\
If 
$$u=\frac{Q(q)}{Q(q^2)}, v=\frac{Q(q^3)}{Q(q^6)}$$ 
then
\begin{equation}
v^4+u^4-v^3u^3+6v^2u^2-16vu=0 
\end{equation}      
\textbf{Proof.}\\ 
From relations (87),(43) and related theory we get
$$
4\sqrt{vu}+\sqrt{(v^2-4)(u^2-4)}=\sqrt{(v^2+4)(u^2+4)}
$$     
after some simplification we get the result.

\[
\]

\centerline{\bf References}\vskip .2in

\noindent

[1]: M.Abramowitz and I.A.Stegun, 'Handbook of Mathematical Functions'. Dover Publications, New York. 1972.

[2]: C. Adiga, T. Kim. 'On a Continued Fraction of Ramanujan'. Tamsui Oxford Journal of Mathematical Sciences 19(1) (2003) 55-56 Alethia University.

[3]: C. Adiga, T. Kim, M.S. Naika and H.S. Madhusudhan 'On Ramanujan`s Cubic Continued Fraction and Explicit Evaluations of Theta-Functions'. arXiv:math/0502323v1 [math.NT] 15 Feb 2005.

[4]: G.E.Andrews, 'Number Theory'. Dover Publications, New York. 1994.
 
[5]: G.E.Andrews, Amer. Math. Monthly, 86, 89-108(1979).

[6]: B.C.Berndt, 'Ramanujan`s Notebooks Part I'. Springer Verlang, New York (1985)

[7]: B.C.Berndt, 'Ramanujan`s Notebooks Part II'. Springer Verlang, New York (1989)

[8]: B.C.Berndt, 'Ramanujan`s Notebooks Part III'. Springer Verlang, New York (1991)

[9]: Bruce C. Berndt, Heng Huat Chan and Liang-Cheng Zhang, 'Ramanujan`s class invariants and cubic continued fraction'. Acta Arithmetica LXXIII.1 (1995).   

[10]: Heng Huat Chan, 'On Ramanujans Cubic Continued Fraction'. Acta Arithmetica. 73 (1995), 343-355.     

[11]: Heng Huat Chan, Sen-Shan Huang, 'On the Ramanujan-Gollnitz-Gordon Continued Fraction'. The Ramanujan Journal 1, 75-90 (1997). 

[12]: I.S. Gradshteyn and I.M. Ryzhik, 'Table of Integrals, Series and Products'. Academic Press (1980).

[13]: Megadahalli Sidda Naika Mahadeva Naika, Mugur Chinna Swamy Maheshkumar, Kurady Sushan Bairy. 'General Formulas for Explicit Evaluations of Ramanuajn`s Cubic Continued Fraction'. Kyungpook Math. J. 49 (2009), 435-450.  

[14]: L. Lorentzen and H. Waadeland, Continued Fractions with Applications. Elsevier Science Publishers B.V., North Holland (1992).  

[15]: S.H.Son, 'Some integrals of theta functions in Ramanujan's lost notebook'. Proc. Canad. No. Thy Assoc. No.5 (R.Gupta and K.S.Williams, eds.), Amer. Math. Soc., Providence.

[16]: H.S. Wall, 'Analytic Theory of Continued Fractions'. Chelsea Publishing Company, Bronx, N.Y. 1948.  

[17]:E.T.Whittaker and G.N.Watson, 'A course on Modern Analysis'. Cambridge U.P. (1927)

[18]:I.J. Zucker, 'The summation of series of hyperbolic  functions'. SIAM J. Math. Ana.10.192. (1979)

[19]: Nikos Bagis and M.L. Glasser, 'Integrals related with Rogers Ramanujan continued fraction and q-products'. arXiv:0904.1641. (2009) 

[20]: W. Duke. 'Continued fractions and Modular functions'. Bull. Amer. Math. Soc. (N.S.), 42 (2005), 137-162. 

\end{document}